\documentclass[12pt,twoside]{amsart}
\usepackage{mathrsfs}
\usepackage{mathtools}
\usepackage{amssymb}
\usepackage{verbatim}
\usepackage{t1enc}
\usepackage{amsmath}
\usepackage{comment}
\usepackage{amsthm,thmtools,xcolor}
\usepackage[colorlinks,linkcolor=blue,citecolor=blue, hypertexnames=false,pdfstartview=FitH]
{hyperref}
\usepackage{backref}
\usepackage{bm}
\usepackage{a4wide}
\usepackage{babel}
\usepackage[latin1]{inputenc}
\usepackage[T1]{fontenc}
\usepackage{times}
\usepackage{hyperref}
\usepackage{amssymb,latexsym}
\usepackage{enumerate}

\usepackage{etoolbox}
\makeatletter
\patchcmd\maketitle
{\uppercasenonmath\shorttitle}
{}
{}{}
\patchcmd\maketitle
{\@nx\MakeUppercase{\the\toks@}}
{\the\toks@}
{}
{}{}
\patchcmd\@settitle
{\uppercasenonmath\@title}
{}
{}{}
\patchcmd\@setauthors
{\MakeUppercase{\authors}}
{\authors}
{}{}

\usepackage{enumerate}

\newcommand{\C}{\ensuremath{\mathbb{C}}}

\newcommand{\D}{\ensuremath{\mathbb{D}}}

\makeatletter
\newcommand{\sumprime}{\if@display\sideset{}{'}\sum%
	\else\sum'\fi}
\makeatother

\begin{document}

	\begin{sloppypar}	
		
		\numberwithin{equation}{section}
		
		% define theorem environments
		\newtheorem{theorem}{Theorem}[section]
		\newtheorem{proposition}[theorem]{Proposition}
		\newtheorem{conjecture}[theorem]{Conjecture}
		\def\theconjecture{\unskip}
		\newtheorem{corollary}[theorem]{Corollary}
		\newtheorem{lemma}[theorem]{Lemma}
		\newtheorem{observation}[theorem]{Observation}
		\newtheorem{definition}[theorem]{Definition}
		
		\newtheorem{remark}{Remark}
		\def\theremark{\unskip}
		\newtheorem{kl}{Key Lemma}
		\def\thekl{\unskip}
		\newtheorem{question}{Question}
		\def\thequestion{\unskip}
		\newtheorem{example}{Example}
		\def\theexample{\unskip}
		\newtheorem{problem}{Problem}
		
		\bigskip
		
		\address{DEPARTMENT OF MATHEMATICAL SCIENCES, NORWEGIAN UNIVERSITY OF SCIENCE AND TECHNOLOGY, NO-7491 TRONDHEIM, NORWAY} 
		
		\email{xu.wang@ntnu.no}
		
		\address{SCHOOL OF MATHEMATICAL SCIENCES, FUDAN UNIVERSITY, SHANGHAI, 20043, CHINA}
		\email{23110180045@m.fudan.edu.cn}

		\title{Hyperconvexity of pseudoconvex domains with continuous boundaries}
		
		\author{Xu Wang and Chenwei Yang}
		\date{\today}

		\begin{abstract} This is a reading note on a disc formula proof of hyperconvexity of bounded pseudoconvex domains in $\mathbb C^n$ with continuous 
			boundaries. 
		\end{abstract}

		\maketitle

		\section{Introduction}  
		
A domain $\Omega\subset\C^n$ is said to be pseudoconvex if there exists a continuous plurisubharmonic (psh) function $\psi$ on $\Omega$ such that $\{\psi<t\}\Subset\Omega$ for every $t\in \mathbb R$ ( we call such $\psi$ an exhaustion). If we can choose $\psi$ such that $\psi=-\log(-\rho)$ for some negative continuous psh $\rho$ on $\Omega$, then we call $\Omega$ a hyperconvex domain. Note that $\{\rho<c\}\Subset\Omega$ for every $c<0$ , we call such $\rho$ a negative exhaustion.  Thus a domain is hyperconvex if it has a continuous
		negative psh exhaustion.
		
		\medskip
		
It is known that bounded pseudoconvex domains are hyperconvex under certain 
smooth boundary assumptions. Diederich--Forn{\ae}ss \cite{DF} proved the
$C^2$ case, then Kerzman--Rosay \cite{KR} the $C^1$ case, Demailly \cite{Demailly} the
Lipschitz case, Avelin--Hed--Persson \cite{AHP}  the log-Lipschitz case, Wang \cite{WangHyper} the radial log-Lipschitz case and finally the H\"older case by Chen \cite{Chen}. 

\medskip

The following continuous boundary case is asked in \cite[Problem~1]{Chen} and  \cite[Problem~2]{Blocki}. 

\begin{theorem}\label{th:hyperconvex} Bounded pseudoconvex domains with continuous boundaries are hyperconvex.	
\end{theorem}

\medskip

\noindent
\textbf{On the use of AI.} When Wang started to use AI in August 2026, he asked this hyperconvexity question to GPT-5.6 Sol and obtained a proof of the above theorem using a precise Berndtsson/Hamano--Maitani--Yamaguchi type curvature formula (see the appendix A,B,C in \cite{Wang-AI}) for variation of bounded simply connected domains in $\mathbb C$. The key in this curvature formula proof is a psh property for the variation of the nome function (see Theorem 7.1 in \cite{Wang-AI}). After that GPT-5.6 found a direct proof that psh property using an extra double cover trick based on the Joukowski map \eqref{eq:Joukowski} (see sections 7.4-7.5 in \cite{Wang-AI}). Since all bounded simply connected domains in $\mathbb C$ are conformal to the disc, one may ask whether there are other "variation of discs" proofs. Then Yang found a direct Poisson disc formula proof (see Proposition 3.2 in \cite{Wang-AI}) using GPT-5.6. Wang thought that using Green disc formula would simplify the proof. Eventually GPT-6 Astra found this simplified proof (see Proposition 3.4 in \cite{Wang-AI}). One interesting thing is, the Joukowski map \eqref{eq:Joukowski} is one of the key ingredients in almost all the proofs (except the first curvature formula proof). This paper is just a reading note for Proposition 3.4 in \cite{Wang-AI}.

\medskip

\noindent
\textbf{Acknowledgements.} It is a pleasure to thank Bo-Yong Chen for sharing this question and pointing out that \eqref{eq:claim} and a recent result of Yuanpu Xiong together imply that the hyperconvexity index of a pseudoconvex domain with  H\"older continuous boundary is positive. In particular, \cite[Corollary 1.8]{Chen17} implies that for fixed $z_0$ in a H\"older bounded pseudoconvex domain, the Bergman distance $ d(z,z_0)$ satisfies that $ d(z,z_0) \geq C |\log \delta(z)|/ \log |\log \delta(z)|$ for $z$ sufficiently close to the boundary, where $\delta$ denotes the boundary distance. We also thank W\l odzimierz Zwonek for his clarifications regarding the proof of Theorem \ref{thm:green-disc-formula} and Xu Xing for explaining the relation between local hyperconvexity index and global hyperconvexity index.

\section{A Green disc formula proof} Let $D$ be a bounded domain. Fix $a\in D$, a canonical negative psh function on $D$ is the pluricomplex Green function with pole at $a$ defined by
$$
G_D(z,a):=\sup \left\lbrace \psi(z): \psi \, \text{is negative psh on $D$ and $\nu_a (\psi) \geq 1$}\right\rbrace,
$$
where
$$
\nu_a (\psi):=\liminf_{z\to a} \frac{\psi(z)}{\log|z-a|}
$$
denotes the Lelong number for $\psi$ at $a$. Denote by $\mathcal O(\bar{\mathbb D}, D)$ the space of holomorphic mappings $f$ on a neighborhood of the closure of the unit disc $\mathbb D$ with $f(\bar{\mathbb D})\subset D$. For $f(\alpha)=a$ we use $m_f(\alpha)$ to denote the vanishing order of $f-a$ at $\alpha$. Then we have the following Poletsky type disc formula for $G_D$ (see \cite{Poletsky1991,Poletsky1993} for the Poletsky theory, \cite{Edigarian97} for the first complete proof,  
for the manifold version see \cite{Rosay, LS2003}, a singular space
formulation is \cite[Theorem~1.1]{DDF2012},  see also \cite{Edwards1966} for an early weak Legendre duality version of the Poletsky theory).

\begin{theorem}[see Theorem 1 in  \cite{Edigarian97}]\label{thm:green-disc-formula}
	For $z\ne a$ in a bounded domain $D\subset\C^n$, we have
	\begin{equation}\label{eq:green-disc-formula}
		G_D(z,a)=
		\inf_{\substack{f\in\mathcal O(\bar{\mathbb D}, D)\\ f(0)=z}}
		\sum_{\substack{\alpha\in\D\\ f(\alpha)=a}}
		m_f(\alpha)\log|\alpha|,
	\end{equation}
	(it is always a finite sum since there are only finite $\alpha \in \mathbb D$ with  $f(\alpha) =a$).
\end{theorem}

To prove Theorem \ref{th:hyperconvex}, we also need the following continuity result for pseudoconvex domains.

\begin{lemma}\label{lem:disc-continuity}
	Let $V$ be a bounded open set in $\mathbb C$ and let
	$X\subset\C^n$ be pseudoconvex. Assume that
	$(s,\zeta)\mapsto f_s(\zeta)$ is continuous on
	$[0,1]\times\overline V$ and each $f_s$ is holomorphic from $V$ to $\mathbb C^n$.
	If
	$$
	f_0(\overline V)\subset X \ \ \ \text{and}  \ \ \ 
	K:=\bigcup_{0\le s\le1}f_s(\partial V)\Subset X,
	$$
	then $f_s(\overline V)\subset X$ for every $s\in[0,1]$.
\end{lemma}

\begin{proof}
	Choose a continuous plurisubharmonic exhaustion $\psi$ of $X$ and put $$I:=\{s\in [0,1]:f_s(\overline V)\subset X\}.$$  Then $I$ contains $0$ and is
	relatively open by compactness of $\overline V$ and continuity of $f_s(\zeta)$.  For $s\in I$, the
	maximum principle applied to $\psi\circ f_s$ gives
	\[
	f_s(\overline V)\subset \left\lbrace\psi\le \max_K\psi\right\rbrace\Subset X,
	\]
	which implies that $I$ is also relatively closed.  Hence
	$I=[0,1]$.
\end{proof}

\begin{proof}[Proof of Theorem \ref{th:hyperconvex}] Let $\Omega$ be a bounded pseudoconvex domain. Fix $a\in \Omega$, it is known that (see the proof of Theorem 1.4.6 in \cite{BlockiNotes}) $\Omega$ is hyperconvex if and only if 
	\begin{equation}\label{eq:barrier}
		\lim_{z\to z_0} G_\Omega(z, a)=0
	\end{equation}  
	for every $z_0\in \partial \Omega$. The locality theorem of Kerzman--Rosay (see \cite[{Proposition 1.1}]{KR})  further implies that $\Omega$ is hyperconvex if for every $z_0\in \partial\Omega$, 
	\begin{equation}\label{eq:barrier1}
		\lim_{z\to \partial\Omega \cap Q} G_{\Omega\cap Q}(z, a)=0 \ \ \left( \text{i.e. $\lim_{z\to w} G_{\Omega\cap Q}(z, a)=0$ for every $w\in \partial\Omega \cap Q$}\right),
	\end{equation} 
	for some convex open neighborhood $Q$ of $z_0$ and some $a\in \Omega\cap Q$. To prove \eqref{eq:barrier1}, we make a holomorphic affine change of coordinate such that $z_0=0$ and
		\begin{equation}\label{eq:Q}
		\Omega\cap 2Q=\left\lbrace z=(z',z_n)\in 2Q: {\rm Im}\, z_n< g(z', {\rm Re}\, z_n) \right\rbrace,
	\end{equation} 
	where, for some small constant $0<r<1/2$, $Q$ is defined by
	$$
	Q:=\{(z', z_n)\in \mathbb C^{n-1} \times\mathbb C: |z'|<r,\, |{\rm Re}\,z_n|<r, \,|{\rm Im}\,z_n|<1\},
	$$ 
	$g(z', {\rm Re}\, z_n) $ is continuous for $|z'|\leq 2r$, $|{\rm Re}\, z_n|\leq 2r$, $g(0)=0$ and 
		\begin{equation}\label{eq:2Q}
\max_{|z'|\leq 2r, \ |{\rm Re}\, z_n| \leq 2r}|g(z', {\rm Re}\, z_n)| <1.
\end{equation} 	
Put
\begin{equation}\label{eq:D}
	D:=\Omega\cap Q=\left\lbrace z\in Q: {\rm Im}\, z_n<g(z', {\rm Re}\, z_n) \right\rbrace.
\end{equation}  
We know that \eqref{eq:barrier1} holds if we can prove 
		\begin{equation}\label{eq:claim}
0<1-e^{G_{D}(z,a)} <\frac 8{\delta(a)} \sqrt{g(z', {\rm Re}\, z_n)-{\rm Im}\, z_n}, \ \ \ z\in D,
\end{equation} 
where $\delta(a):=\inf_{z\in \partial D} |z-a|$ denotes the boundary distance function. 
	To prove \eqref{eq:claim}, we apply Theorem \ref{thm:green-disc-formula}. For $f\in \mathcal O(\bar{\mathbb D}, D)$ with $f(0)=z$, put
\begin{equation}\label{eq:green-blaschke-product}
	\mathcal B_{f, a}(\zeta):=c \prod_{f(\alpha)=a, \, \alpha\in \mathbb D}
	\left(\frac{\zeta-\alpha}
	{1-\overline{\alpha}\zeta}\right)^{m_f(\alpha)}, \qquad \zeta \in \mathbb C,
\end{equation}
where the constant $|c|=1$ is chosen such that $\mathcal B_{f, a}(0)=\prod_{j=1}^s|\alpha_j|^{m_j}>0$. Since $G_D<0$, we always have $0<1-e^{G_{D}(z,a)} $. Thus Theorem \ref{thm:green-disc-formula} implies that \eqref{eq:claim} is equivalent to
\begin{equation}\label{eq:green-blaschke-product1}
	1- \mathcal B_{f, a}(0) < \frac 8{\delta(a)}  \sqrt{g(z', {\rm Re}\, z_n)-{\rm Im}\, z_n}, \qquad z \in D,
\end{equation}
  for all $f\in \mathcal O(\bar{\mathbb D}, D)$ with $f(0)=z$. \textcolor{blue}{The main idea is to introduce the following deformation
\begin{equation}\label{eq:green-blaschke-product2}
	f_s(\zeta):=f(\zeta)+ \left( 0', -i s \varepsilon^2 B(0) J(B(\zeta))\right), \ \ \varepsilon:=\frac{\delta(a)}{4}<\frac12,
\end{equation}
to $f\in \mathcal O(\bar{\mathbb D}, D)$ (with $f(0)=z$), where $0'$ denotes the origin in $\mathbb C^{n-1}$, $B$ is short for $\mathcal B_{f, a}$ and 
\begin{equation}\label{eq:Joukowski}
J(\lambda):=\frac{2-\lambda-\lambda^{-1}}{4}
\end{equation}
denotes the Joukowski map}. We shall prove \eqref{eq:green-blaschke-product1} by applying Lemma \ref{lem:disc-continuity} to $X:=\Omega\cap 2Q$ and 
$
V:=\left\lbrace \zeta\in \mathbb D: |B(\zeta)|>\varepsilon B(0)\right\rbrace.
$ Let us verify the assumptions in Lemma \ref{lem:disc-continuity}: $f_0(\overline V)\subset X $ holds since $D\subset X$; to verify the second, we use 
$\partial V \subset \{|\zeta|=1\}\cup\left\lbrace |B(\zeta)|=\varepsilon B(0)\right\rbrace.
$ 
For the first $|\zeta|=1$ part we have $|B(\zeta)|=1$. Since $J$ sends $|\lambda|=1$ to $[0,1]$, \eqref{eq:Q} and \eqref{eq:2Q} give
$$
\bigcup_{0\le s\le1} f_s(\partial\mathbb D)\in f(\partial\mathbb D) + \{(0', -i t\varepsilon^2): t\in [0,1]\} \Subset \Omega\cap 2Q.
$$
To handle the second part of $\partial V$, we apply the maximum principle to $(f(\zeta)-a)/B(\zeta)$:
$$
|f(\zeta)-a| \leq |B(\zeta)| \sup_{|\xi|= 1} \big |\frac{f(\xi)-a}{B(\xi)}| = |B(\zeta)| \sup_{|\xi|= 1} \big |f(\xi)-a|   \leq 3 |B(\zeta)|, \ \ \forall \ |\zeta|\leq 1.
$$
Thus for $\zeta$ with $|B(\zeta)|=\varepsilon B(0)  \leq \varepsilon =\delta(a)/4 $, we have  $|f(\zeta)-a| \leq 3\delta(a)/4 $ and
$$
|f_s(\zeta)-a| \leq |f(\zeta)-a|+ \varepsilon |B(\zeta) J(B(\zeta))| \leq \frac{3 \delta(a) }{4} +  \frac{\delta(a)}{4} \frac{9}{16} <\delta(a), \ \ \ \ \forall \ s\in [0,1].
$$
Thus we have $\bigcup_{0\le s\le1}f_s(\partial V)\Subset\Omega \cap 2Q$ and Lemma \ref{lem:disc-continuity} for $X=\Omega\cap 2Q$ gives $f_s(\bar V) \subset X$ for all $s\in [0,1]$. In particular, 
$$
f_1(0) = \left( z', z_n + \frac{i \varepsilon^2 (1-B(0))^2}{4}\right) \in \Omega \cap 2Q
$$
gives
$
{\rm Im}\, z_n + \frac{\varepsilon^2 (1-B(0))^2}{4} < g(z', {\rm Re}\, z_n),
$
from which \eqref{eq:green-blaschke-product1} follows.
\end{proof}

\end{sloppypar}


\begin{thebibliography}{99}
	

	
	\bibitem{AHP}
	B.~Avelin, L.~Hed and H.~Persson,
	\emph{A note on the hyperconvexity of pseudoconvex domains beyond
		Lipschitz regularity}, Potential Anal. \textbf{43} (2015), no.~3,
	531--545.
	\href{https://doi.org/10.1007/s11118-015-9486-1}
	{doi:10.1007/s11118-015-9486-1}.

	\bibitem{BlockiNotes}
	Z.~B\l ocki, \emph{The complex Monge--Amp\`ere operator in
		pluripotential theory}, lecture notes, Jagiellonian University,
	1998. \href{https://gamma.im.uj.edu.pl/~blocki/publ/ln/wykl.pdf}
	{Author's online notes}.
	
	\bibitem{Blocki}
	Z.~B\l ocki, \emph{Bergman kernel and pluripotential theory}, in
	Analysis, Complex Geometry, and Mathematical Physics, Contemp. Math.
	\textbf{644}, Amer. Math. Soc., 2015, 1--10.
	\href{https://doi.org/10.1090/conm/644/12785}{doi:10.1090/conm/644/12785}.
	
	
	\bibitem{Chen17} B.-Y. Chen, {\it Bergman kernel and hyperconvexity index}, Anal. PDE {\bf 10} (2017), 1429--1454.
	

	
	\bibitem{Chen}
	B.-Y.~Chen, \emph{Every bounded pseudoconvex domain with H\"older
		boundary is hyperconvex}, Bull. Lond. Math. Soc. \textbf{53} (2021),
	1009--1015.
	\href{https://doi.org/10.1112/blms.12478}{doi:10.1112/blms.12478}.
	
	\bibitem{Demailly}
	J.-P.~Demailly, \emph{Mesures de Monge--Amp\`ere et mesures
		pluriharmoniques}, Math. Z. \textbf{194} (1987), 519--564.
	\href{https://doi.org/10.1007/BF01161920}{doi:10.1007/BF01161920}.
	
	
	\bibitem{DF}
	K.~Diederich and J.~E.~Forn\ae ss, \emph{Pseudoconvex domains:
		bounded strictly plurisubharmonic exhaustion functions}, Invent. Math.
	\textbf{39} (1977), 129--141.
	
	\bibitem{DDF2012}
	B.~Drinovec Drnov\v{s}ek and F.~Forstneri\v{c},
	\emph{The Poletsky--Rosay theorem on singular complex spaces},
	Indiana Univ. Math. J. \textbf{61} (2012), no.~4, 1407--1423.
	\href{https://users.fmf.uni-lj.si/forstneric/papers/2012Indiana.pdf}
	{Author's copy}; see Theorem~1.1.
	
	\bibitem{Edwards1966}
	D.~A.~Edwards, \emph{Choquet boundary theory for certain spaces of
		lower semicontinuous functions}, in \emph{Function Algebras}
	(Proc. Internat. Sympos., Tulane University, 1965),
	F.~T.~Birtel (ed.), Scott, Foresman and Co., Chicago, 1966,
	300--309.
	
	\bibitem{Edigarian97} A. Edigarian, {\it On definitions of the pluricomplex Green function}, Ann. Polon. Math. {\bf 67} (1997), 233--246.
	
	
	\bibitem{Wang-AI}
	GPT-5.6 Sol and GPT-6 Astra,
	\emph{The elliptic nome, condenser capacity, and hyperconvexity}, see the  \href{https://xuwa.folk.ntnu.no/res/nome.pdf}{link}.

	\bibitem{KR}
	N.~Kerzman and J.-P.~Rosay,
	\emph{Fonctions plurisousharmoniques d'exhaustion born\'ees et domaines
		taut}, Math. Ann. \textbf{257} (1981), 171--184.
	\href{https://doi.org/10.1007/BF01458282}{doi:10.1007/BF01458282}.
	
	\bibitem{LS2003}
	F.~L\'arusson and R.~Sigurdsson,
	\emph{Plurisubharmonicity of envelopes of disc functionals on
		manifolds}, J. reine angew. Math. \textbf{555} (2003), 27--38.
	\href{https://doi.org/10.1515/crll.2003.013}{doi:10.1515/crll.2003.013}.
	

	\bibitem{Poletsky1991}
	E.~A.~Poletsky, \emph{Plurisubharmonic functions as solutions of
		variational problems}, in \emph{Several Complex Variables and Complex
		Geometry, Part~1} (Santa Cruz, CA, 1989), Proc. Sympos. Pure Math.
	\textbf{52}, Part~1, Amer. Math. Soc., Providence, RI, 1991, 163--171.
	\href{https://doi.org/10.1090/pspum/052.1/1128523}
	{doi:10.1090/pspum/052.1/1128523}.
	
	\bibitem{Poletsky1993}
	E.~A.~Poletsky, \emph{Holomorphic currents}, Indiana Univ. Math. J.
	\textbf{42} (1993), no.~1, 85--144.
	\href{https://doi.org/10.1512/iumj.1993.42.42006}
	{Journal record}.
	
	\bibitem{Rosay}
	J.-P.~Rosay, \emph{Poletsky theory of disks on holomorphic
		manifolds}, Indiana Univ. Math. J. \textbf{52} (2003), no.~1,
	157--169.
	\href{https://doi.org/10.1512/iumj.2003.52.2170}
	{doi:10.1512/iumj.2003.52.2170}.
	
	\bibitem{WangHyper}
	X. Wang,
	\emph{Hyperconvexity of non-smooth pseudoconvex domains},
	Ann. Polon. Math. \textbf{111} (2014), 1--11.
	


	
\end{thebibliography}
\end{document}